\documentclass[12pt]{amsart}
\usepackage{fullpage}
\usepackage{enumitem}
\usepackage{comment}
\newtheorem{theorem}{Theorem}
\newtheorem{lemma}[theorem]{Lemma}
\newtheorem{proposition}[theorem]{Proposition}

\newtheorem{conjecture}[theorem]{Conjecture}
\newtheorem{claim}{Claim}[theorem]

\newtheorem*{claim*}{Claim}

\theoremstyle{definition}

\newtheorem{remark}[theorem]{Remark}
\newtheorem{definition}[theorem]{Definition}
\newtheorem{example}[theorem]{Example}

\newcommand{\cA}{\mathcal{A}}

\newcommand{\cF}{\mathcal{F}}

\newcommand{\cP}{\mathcal{P}}

\newcommand{\cT}{\mathcal{T}}

\newcommand{\floor}[1]{\left\lfloor #1 \right\rfloor}
\newcommand{\ceiling}[1]{\left\lceil #1 \right\rceil}

\title{Transversal $C_k$-factors in subgraphs of the balanced blow-up of $C_k$}
\author{Beka Ergemlidze$^{1}$, Theodore Molla$^{2}$}

\date{\today}

\begin{document}

\begin{abstract}
  For a subgraph $G$ of the blow-up of a graph $F$,
  we let $\delta^*(G)$ be the smallest minimum degree over all of the 
  bipartite subgraphs of $G$ induced by pairs of parts 
  that correspond to edges of $F$.
  In
  [Triangle-factors in a balanced blown-up triangle. Discrete Mathematics, 2000],
  Johansson proved that if $G$ is a spanning subgraph of the blow-up of $C_3$ with
  parts of size $n$ and $\delta^*(G) \ge \frac{2}{3}n + \sqrt{n}$, 
  then $G$ contains $n$ vertex-disjoint triangles,
  and presented the following conjecture of H\"aggkvist:
  If $G$ is a spanning subgraph of the blow-up of $C_k$ with parts
  of size $n$ and $\delta^*(G) \ge (1 + 1/k)n/2 + 1$, then
  $G$ contains $n$ vertex disjoint copies of $C_k$ such that each $C_k$ intersects
  each of the $k$ parts exactly once.
  The degree condition of this conjecture is tight when $k=3$ and 
  cannot be strengthened by more than one when $k \ge 4$.
  A similar conjecture was also made by Fischer in
  [Variants of the Hajnal-Szemer\'edi Theorem.  Journal of Graph Theory, 1999]
  and the triangle case was proved for large $n$ by Magyar and Martin in 
  [Tripartite version of the Corr\'adi-Hajnal Theorem.  Discrete Mathematics, 2002].
  
 In this paper, we prove this Conjecture asymptotically.
  We also pose a conjecture which generalizes this result
  by allowing the minimum degree conditions on the nonempty bipartite subgraphs 
  induced by pairs of parts to vary. 
  Our second result supports this new conjecture by proving the triangle case.
  This result generalizes Johannson's result asymptotically.  
\end{abstract}

\maketitle

\noindent\footnotetext[1]{Department of Mathematics and Statistics, University of South Florida {\tt ergemlidze@usf.edu}.}  

\noindent\footnotetext[2]{Department of Mathematics and Statistics, University of South Florida {\tt molla@usf.edu}. 
Research supported in part by NSF Grant DMS~1800761.}

\section{Introduction}

For a graph $F$ on $[k] := \{1, \dotsc, k\}$, we say that
$B$ is the \textit{$n$-blow-up of $F$} 
if there exists an ordered partition $(V_1, \dotsc, V_k)$ of $V(B)$
such that $|V_1| = \dotsm = |V_k| = n$ and we have that
$uu' \in E(B)$ if and only if $u \in V_i$ and $u' \in V_{j}$ for some 
$ij \in E(F)$.
For $G$ a spanning subgraph of $B$, 
we call the sequence $V_1, \dotsc, V_k$ the \textit{parts} of $G$ and we define 
\begin{equation*}
  \delta^*_F(G) := \min_{ij \in E(F)} \delta(G[V_i, V_{j}])
\end{equation*}
where $G[V_i, V_j]$ is the bipartite subgraph of $G$ induced by the parts 
$V_i$ and $V_j$. 
We often drop the subscript $F$ when it is clear from the context.
For a graph $H$, we call $\cT$ an $H$-tiling of $G$
if $\cT$ consists of vertex disjoint copies of $H$ in $G$.
We say that $\cT$ \textit{covers} 
$V(\cT) := \bigcup\{ V(H') : H' \in \cT \}$ and 
say that $\cT$ is \textit{perfect} or an
\textit{$H$-factor} if it covers every vertex of $G$.
Call a subset of $V(G)$ or a subgraph of $G$ 
a \textit{transversal} if it intersects each
part in exactly one vertex and 
a \textit{partial transversal} if it intersects each
part in at most one vertex. 
An $H$-tiling is a \textit{transversal $H$-tiling} 
if each copy of $H$ in $\cT$ is a transversal.
We call a perfect transversal $H$-tiling a \textit{transversal $H$-factor}.

Fischer \cite{fischer1999} conjectured 
the following multipartite version of the Hajnal-Szemer\'edi Theorem:
If $G$ is the $n$-blow-up of $K_k$, and
$\delta^*(G) \ge \left(1 - \frac{1}{k}\right)n$, then $G$ has
a $K_k$-factor.
In the same paper, Fischer proved that, when $k\in \{3, 4\}$, such
a graph $G$ contains a $K_k$-tiling of size at least 
$n-C$, where $C$ is a constant that depends only on $k$.
Johansson \cite{johansson2000triangle} proved that,
for every $n$, if $G$ is a spanning subgraph of the $n$-blow-up of $K_3$
and $\delta^*(G) \ge 2n/3 + \sqrt{n}$, then $G$ contains a $K_3$-factor,
so Johansson proved the triangle case of the conjecture asymptotically.
Later,
Lo \& M\"arkstrom \cite{lo2013} and, independently, Keevash \& Mycroft \cite{keevash2014}
proved the conjecture asymptotically for every $k \ge 4$.
The following theorem,
which was proved for $k=3$ by Magyar \& Martin \cite{magyar2002},
for $k=4$ by Martin \& Szemer\'edi \cite{martin2008},
and for $k \ge 5$ by Keevash \& Mycroft
\cite{keevash2014},
shows that Fischer's original conjecture was nearly true for $n$ sufficiently
large.
(Keevash \& Mycroft 
actually proved more, see Theorem 1.1 in \cite{keevash2014} for details.)
\begin{theorem}\label{thm:fischer}
  For every $k$ there exists $n_0 := n_0(k)$ such that whenever 
  $n \ge n_0$ the following holds
  for every spanning subgraph $G$ of the $n$-blow-up of $K_k$ where
  \begin{equation*}
    \delta^*(G) \ge \left(1 - \frac{1}{k}\right)n. 
  \end{equation*}
  The graph $G$ does not contains a $K_k$-factor if and only if 
  both $n$ and $k$ are odd, $k$ divides $n$ and $G$ is isomorphic to a specific
  spanning subgraph $\Gamma_{n,k}$ of the $n$-blow-up of $K_k$ 
  where $\delta^*(\Gamma_{n,k}) = \left(1 - \frac{1}{k}\right)n$. 
\end{theorem}

The following conjecture of  H\"aggkvist, which appeared in \cite{johansson2000triangle}, 
can be seen as a different generalization of the $k=3$ 
case of Theorem~\ref{thm:fischer}.
Independently, Fischer made a similar conjecture in \cite{fischer1999}.
\begin{conjecture}\label{conj:main}
  For every $k \ge 3$, if $G$ is a spanning subgraph of the $n$-blow-up of $C_k$ and 
  \begin{equation}\label{eq:conj_degree_condition}
    \delta^*(G) \ge \left(1 + \frac{1}{k}\right)\frac{n}{2} + 1, 
  \end{equation}
  then $G$ has a transversal $C_k$-factor. 
\end{conjecture}

Our first result establishes an asymptotic version of Conjecture~\ref{conj:main}.
\begin{theorem}\label{thm:asymp}
  For every $\varepsilon > 0$ and positive integer $k \ge 4$ 
  there exists $n_0 := n_0(k,\varepsilon)$ such that 
  for every $n \ge n_0$ the following holds.
  If $G$ is a spanning subgraph of the $n$-blow up of $C_k$ and 
  \begin{equation}\label{eq:asymp_degree_condition}
    \delta^*(G) \ge \left(1 + \frac{1}{k} + \varepsilon\right)\frac{n}{2}, 
  \end{equation}
  then $G$ has a transversal $C_k$-factor.
\end{theorem}

Note that Theorem~\ref{thm:fischer} shows that Conjecture~\ref{conj:main} is
tight when $k = 3$.
The following example from \cite{johansson2000triangle} 
shows that, for $k \ge 4$,
the minimum degree condition \eqref{eq:conj_degree_condition} 
in Conjecture~\ref{conj:main} 
cannot be decreased by more than $1$. 
Call $Z \subseteq V(G)$ a \textit{transversal $C_k$-cover} if
every transversal $C_k$ in $G$ intersects $Z$ and let the 
\textit{transversal $C_k$-cover number} of $G$ be the order
of a smallest transversal $C_k$-cover.
This example relies on the observation that the
maximum size of a transversal $C_k$-tiling is bounded above
by the transversal $C_k$-cover number.
(Note that we always view arithmetic on elements of $[k] := \{1, \dotsc, k\}$ modulo $k$.)
\begin{example}\label{exp:conj}
  For $k \ge 3$ and $m \ge 1$, let $n := 2km$ and
  $V_1, \dotsc, V_k$ be disjoint sets each of size $n$.
  For $i \in [k-1]$, let $\{U_i, W_i, Z_i\}$ be a partition of $V_i$
  such that $|U_i| = (k-1)m$, $|W_i| = (k-1)m$ and $|Z_i| = 2m$, and let $\{U_k, W_k, Z_k\}$ be a partition of $V_k$ such that $|U_k| = (k-1)m$, $|W_k| = (k-1)m+1$ and $|Z_k| = 2m-1$.
  Let $G$ be the spanning subgraph of the $n$-blow-up of $C_k$ with parts 
  $V_1,\dotsc,V_k$ where $E(G)$ consists of the union of the 
  edges in the following graphs:
  \begin{itemize}
    \item the complete bipartite graphs with parts $Z_i, V_{i-1}$ and $Z_i, V_{i+1}$ for each $i \in [k]$,
    \item the complete bipartite graphs with parts $U_i,U_{i+1}$ and $W_i,W_{i+1}$ for each $i \in [k-1]$,
    \item the complete bipartite graphs with parts $U_k,W_1$ and $W_k,U_1$. 
  \end{itemize}
  Note that 
  $\delta^*(G) = (k+1)m - 1 = \left(1 + \frac{1}{k}\right)\frac{n}{2} - 1$, 
  and that every transversal $C_k$ has at least one vertex in 
  $Z := Z_1 \cup \dotsm \cup Z_k$, i.e.,
  $Z$ is a transversal $C_k$-cover of $G$.
  The fact that $|Z| = 2mk - 1 < n$ then implies that $G$ does not contain a 
  transversal $C_k$-factor. 
\end{example}

\bigskip
We make the following conjecture which is strengthening of Theorem~\ref{thm:asymp}.

\begin{conjecture}\label{conj:strong}
  For every $k \ge 3$ and $\varepsilon > 0$, there exists $n_0 := n_0(k,\varepsilon)$ 
  such that for every $n \ge n_0$ the following holds. 
  Let $G$ be a spanning subgraph of the $n$-blow-up of $C_k$ with parts $V_1, \dotsc, V_k$. If there exist $\delta_1,\delta_2, \ldots, \delta_k\ge n/2$ such that $\delta(G[V_i,V_{i+1}])\geq \delta_i$, for every $i\in[k]$, and 
  \begin{equation*}
   \frac1k \sum_{i\in[k]}\delta_i \ge \left(1 + \frac{1}{k}+\varepsilon\right)\frac{n}{2},
  \end{equation*}
  then $G$ has a transversal $C_k$-factor.
\end{conjecture}

Note that Theorem \ref{thm:asymp} is a special, uniform case of Conjecture \ref{conj:strong}, namely the case of $\delta_1=\delta_2=\ldots =\delta_k$.  Also, note that the condition 
$\delta_1, \dotsc, \delta_k \ge n/2$ is necessary
because a transversal $C_k$-factor in $G$ defines a perfect matching
in $G[V_i, V_{i+1}]$ for every $i \in [k]$ and $n/2$ is the smallest minimum degree
condition necessary to guarantee a perfect matching in a bipartite graph with parts of size $n$.

Our second result shows that that Conjecture \ref{conj:strong} holds for $k=3$. 
\begin{theorem}\label{thm:strongk3}
For every $\varepsilon > 0$ there exists $n_0 := n_0(k,\varepsilon)$ such that 
  for every $n \ge n_0$ the following holds.
   Let $G$ be a spanning subgraph of the $n$-blow-up of a triangle with parts $V_1,V_2,V_3$. If there exist $\delta_1,\delta_2,\delta_3\ge n/2$ such that $\delta(G[V_i,V_{i+1}])\geq \delta_i$, for every $i\in[3]$, and 
  \begin{equation*}
    \frac{\delta_1 + \delta_2 + \delta_3}3
    \ge \left(\frac{2}{3}+\varepsilon\right)n,
  \end{equation*}
  then $G$ has a triangle factor.
\end{theorem}
It is interesting to observe that in the case of $\delta_i\le (1 + \varepsilon)\frac n2$ for some $i\in[k]$, the problem in Conjecture \ref{conj:strong} for $k$ can be reduced to $k-1$. Indeed, let us consider a spanning subgraph $G$ of the $n$-blow-up of $C_{k}$ where $\delta_k=(1+\varepsilon) \frac n2$ (for $i\in[k]$, $\delta_i$ and $V_i$ are defined as in conjecture). Hall's Theorem implies that we can match every $v \in V_{1}$ 
  to a unique $f_v \in V_{k}$ that is adjacent to $v$. 
  Let $G'$ be the graph derived from $G$ by collapsing each edge $vf_v$ into $v$ for each $v\in V_1$, i.e., 
  $G'$ is $G - V_k$ with an edge between $v \in V_1$ and $u \in V_{k-1}$ if and only if $f_v$ is adjacent to $u$ in $G$.
It is easy to see that $G'$ is a spanning subgraph of the $n$-blow-up of $C_{k-1}$ 
  with parts $V_1, \dotsc, V_{k-1}$ such that $\delta(G[V_i,V_{i+1}]) \ge \delta_i$ for each $i\in [k-1]$,
  \begin{equation*}
    \sum_{i \in [k-1]} \delta_i \ge 
    \sum_{i \in [k]} \delta_i - (1+\varepsilon)\frac n2 \ge 
    k\left(1 + \frac{1}{k} + \varepsilon\right)\frac n2 - (1+\varepsilon)\frac n2
    = (k-1)\left(1 + \frac{1}{k-1} + \varepsilon\right) \frac n2, 
  \end{equation*}
 and that any transversal $C_{k-1}$-factor in $G'$ can be extended to a transversal $C_{k}$-factor in $G$.

The above observation and Theorem \ref{thm:strongk3} imply that Conjecture \ref{conj:strong} holds if $\delta_1,\delta_2,\ldots, \delta_k\ge n/2$ and $\delta_i+\delta_j+\delta_l\ge 2n+\varepsilon$ (i.e., in the case when all the excess values of $\delta_i$ compared to $n/2$ are concentrated at $3$ members of $\delta_1,\delta_2,\ldots \delta_k$).

We also note that this observation and 
a straightforward application of the absorbing method of R\"odl, Ruci\'nski, and Szemer\'edi
together imply that one only needs to prove the following weaker conjecture to establish
Conjecture~\ref{conj:strong}.
(This reduction is proved in Lemma~\ref{lem:weakerstrongtostrong} in Section~\ref{sec:absorbing}.) 
\begin{conjecture}\label{conj:weakerstrong}
  For every $k \ge 3$ and $\varepsilon > 0$, there exists $\sigma_0 := \sigma_0(k, \varepsilon)$
  such that for $0 < \sigma < \sigma_0$ there exists $n_0 := n_0(k,\varepsilon, \sigma)$ 
  such that for every $n \ge n_0$ the following holds. 
  Let $G$ be a spanning subgraph of the $n$-blow-up of $C_k$ with parts $V_1, \dotsc, V_k$. If there exist $\delta_1,\delta_2, \ldots, \delta_k\ge (1 + \varepsilon) n/2$ such that $\delta(G[V_i,V_{i+1}])\geq \delta_i$, for every $i\in[k]$, and 
  \begin{equation*}
   \frac1k \sum_{i\in[k]}\delta_i \ge \left(1 + \frac{1}{k}+\varepsilon\right)\frac{n}{2},
  \end{equation*}
  then $G$ has a transversal $C_k$-tiling of size at least $(1-\sigma)n$.
\end{conjecture}

To further support the conjecture, we consider a natural extension of Conjecture \ref{conj:strong} for $k=2$.
  Suppose that $H$ and $H'$ are two balanced bipartite graphs both with the same partite sets $V_1$ and $V_2$
  where $|V_1| = |V_2| = n$.
  If $\frac 12\left(\delta(H) + \delta(H')\right) \ge \frac{3n}{4}$, 
  then there exists $M \subseteq E(H) \cap E(H')$ that
  is simultaneously a perfect matching of both $H$ and $H'$. Indeed, every vertex in $v \in V_1 \cup V_2$ is incident to at least $d_H(v) + d_{H'}(v) - n \ge n/2$ edges
  in $E(H) \cap E(H')$, so Hall's Theorem implies the desired matching exists.
 This together with the reduction argument discussed above, implies the following.
  \begin{remark}\label{remark:C2}
For every $k \ge 3$ and $n$, Conjecture~\ref{conj:strong} holds
  whenever $\delta_1, \dotsc, \delta_k\ge n/2$ and 
  $\delta_i + \delta_j \ge 3n/2$ for distinct $i,j \in [k]$.
\end{remark}
\bigskip

Because of Example~\ref{exp:conj}, 
the condition on the average of the minimum degrees 
in Conjecture~\ref{conj:strong} is asymptotically sharp.
However, it might be possible to weaken the degree condition 
by only placing a lower bound on the average of some
proper subset of the minimum degrees.
For example, in the triangle case, we do not have an example in
which all of the minimum degrees are at least
$n/2$ and the average of only the two largest minimum degrees is at least $2n/3$
that does not have a triangle factor.
Often one tries to find such examples
that either have an independent set which is larger than $n$
or have a triangle cover of size less than $n$, since
either one of these two conditions imply 
that the graph cannot contain $n$ vertex disjoint triangles.
It is a straightforward exercise to show that, under these conditions, 
the independence number must be $n$.
The following theorem proves that the triangle cover number
must be $n$ as well.
\begin{theorem}
\label{thm:trianglecover}
  For every $n \in \mathbb{N}$, the following holds for every spanning subgraph $G$ of the $n$-blow-up
  of $C_3$ with parts $V_1,V_2,V_3$.
  If $\delta_1 \ge \delta_2 \ge \delta_3 \ge n/2$,
  $\delta(G[V_i, V_{i+1}]) \ge \delta_i$ for $i \in [3]$, and
  \begin{equation*}
    \frac{\delta_1 + \delta_2}2 \ge \frac{2n}{3},
  \end{equation*}
  then the triangle cover number of $G$ is $n$. 
  
  Moreover, for every rational $\gamma \in (\frac 34,\frac79] \cup \{\frac23\}$ there are infinitely many $n \in \mathbb{N}$ such that when $\beta = 4/3 - \gamma$ there exists a spanning subgraph $G$ of the $n$-blow-up of $C_3$ with parts $A$, $B$, and $C$ 
  such that $\delta(G[A,B]) \ge \gamma n -1$, $\delta(G[A,C]) \ge \beta n$ and 
  $\delta(G[B,C]) \ge n/2$ that has a triangle cover of order less than $n$.
\end{theorem}

Suppose that for every sufficiently small $\varepsilon > 0$ 
there exists $n_0$ such that when $n \ge n_0$ there
exists a subgraph of the $n$-blow-up of $C_3$ with parts $V_1,V_2,V_3$
that meets the stronger degree conditions 
$\delta_1 \ge \delta_2 \ge \delta_3 \ge (1 + \varepsilon)n/2$ and 
$(\delta_1 + \delta_2)/2 \ge (2/3 + \varepsilon)n$ yet
does not have a triangle factor.
In Section~\ref{sec:absorbing} (Lemma~\ref{lem:absorb} and Proposition~\ref{prop:c3linked}),  
we will show that, under these conditions, 
we can apply the absorbing method.
This would therefore mean that, for some $\sigma > 0$
and for every sufficiently large $n$,
there would exist a subgraph of the $n$-blow-up of $C_3$
that meets the degree conditions of first part of Theorem~\ref{thm:trianglecover}
in which every triangle factor has order at most
$(1 - \sigma)n$, 
but has triangle cover number $n$ and independence number $n$.\footnote{
  It is well-known that this assumption would
  also imply that there would exist a family of examples that 
  meet the degree conditions of first part of Theorem~\ref{thm:trianglecover} and
  do not have a perfect factional triangle tiling.  By the duality
  theorem from linear programming, such a family of examples then
  must have a fractional triangle cover of size less than the size of the parts.}


\subsection{Notation}
For a graph $G$,
$e(G)$ denotes the number of edges in $G$. For $S, T \subseteq V(G)$,  
we let $N_G(S, T) := T \cap \left(\bigcap \{N_G(v) : v \in S \}\right)$
be the \textit{common neighborhood of $S$ in $T$}
and we let $d_G(S, T) := |N_G(S, T)|$.
For $v \in V(G)$, we define $N_G(v, T) := N_G(\{v\}, T)$ 
and $d_G(v, T) = d_G(\{v\}, T)$.
We typically drop the subscript from this notation when it is clear from 
the context.
For a tiling $\cT$, 
we let $U(\cT) := V(G) \setminus V(\cT)$ be the
vertices \textit{uncovered} by $\cT$ and if 
$v \in U(\cT)$ we say that \textit{$v$ is uncovered by $\cT$}.
Similarly, if $e \in E(G)$, and both endpoints of $e$ are uncovered
by $\cT$, we say that \textit{$e$ is uncovered by $\cT$}.

\section{The Absorbing Method}\label{sec:absorbing}
We use a straightforward application of the 
absorbing method of R\"odl, Ruci\'nski, and Szemer\'edi \cite{rodl2006dirac}.
Propositions~\ref{prop:cklinked}  and \ref{prop:c3linked} are essentially 
all that is necessary
to derive appropriate absorbing lemmas in this setting.

\begin{definition}
  For $k \ge 3$, let $G$ be a subgraph of the $n$-blow-up of $C_k$ with parts $V_1, \dotsc, V_k$.
  For vertices $v,v'$ in the same part, we call a $k$-tuple
  of distinct vertices 
  $(v_1,\ldots v_{t})$ a $(v,v', t)$-linking sequence if
  both  $G[\{v,v_1,\dotsc v_t\}]$ and $G[\{v',v_1,\dotsc,v_{t}\}]$ have
  a transversal $C_k$-factor.
  We allow $v = v'$ in this definition.
  We say that $G$ is $(\eta,t)$-linked if, for every $i \in [k]$ and $v,v' \in V_i$,
  the number of $(v,v',t)$-linking sequence is at least $\eta n^{t}$.
\end{definition}

The proof of the following lemma is standard (e.g.\ it is very similar to Lemma 1.1 in \cite{lo2015}), but we include a proof in the appendix for completeness.
\begin{lemma}[The Absorbing Lemma]\label{lem:absorb}
  For $k \ge 3$, $t \ge k-1$, $\eta >0$ and $0 < \sigma \le \frac{0.1 \eta^{k+1}}{(k(t+1))^2 + 1}$,
  there exists $n_0(k, t, \eta, \sigma)$ such that for every $n \ge n_0$ the following holds.
  Suppose that $G$ is a subgraph of the $n$-blow-up of $C_k$ with parts $V_1, \dotsc, V_k$
  that is $(2\eta,t)$-linked. 
  For some $z \le \sigma n$, 
  there exists $A \subseteq V(G)$ where $|A \cap V_i| = z$ for every $i \in [k]$
  such that if $G - A$ has a transversal $C_k$-tiling of size at least
  $n - z - \sigma^2 n$, then $G$ has a transversal $C_k$-factor.
\end{lemma}

Note that the degree condition in the following proposition is weaker than 
the degree condition in Conjecture~\ref{conj:main}. 
\begin{proposition}\label{prop:cklinked}
  For $k \ge 4$ and $\varepsilon > 0$,  
  if $G$ is a subgraph of the $n$-blow-up of $C_k$ and 
  $\delta^*(G) \ge (1 + \varepsilon)n/2$, 
  then $G$ is $(\varepsilon^3/2^{k},k-1)$-linked.
\end{proposition}
\begin{proof}
  Let $V_1, \dotsc, V_k$ be the parts of $G$.
  Without loss of generality we can assume that $v,v' \in V_1$.
  We can construct $(v_2, \dotsc, v_k)$ a $(v,v')$-linking sequence 
  by first selecting $v_2 \in N(\{v, v'\}, V_2)$ and then
  $v_k \in N(\{v, v'\}, V_k)$ each
  in at least $2 \delta^*(G) - n \ge \varepsilon n$ ways.
  Iteratively, for $i$ from $3$ to $k-2$ 
  we can select $v_i \in N(v_{i-1}, V_i)$
  in at least $\delta^*(G) \ge n/2$ ways.
  Finally, we can select 
  $v_{k-1} \in N(\{v_{k-2}, v_k\}, V_{k-1})$ in one at least 
  $2 \delta^*(G) - n \ge \varepsilon n$ ways.
\end{proof}

\begin{proposition}\label{prop:c3linked}
  For every $\varepsilon > 0$ there exist $n_0(\varepsilon)$ such that
  for every $n \ge n_0$ the following holds.
  Let $G$ be a subgraph of the $n$-blow-up of a triangle with parts $V_1,V_2, V_3$. 
  If $\delta_1 \ge \delta_2 \ge \delta_3\geq (1+\varepsilon)n/2$ are such that 
  $\delta(G[V_i,V_{i+1}])\ge\delta_i$ for every $i \in [3]$, and
  \begin{equation*}
    \frac{\delta_1 + \delta_2}{2} \geq \left(\frac23+\varepsilon\right)n,
  \end{equation*}
  then $G$ is $(\varepsilon^3/100,5)$-linked.
\end{proposition}
\begin{proof}
  There are at least 
  $n \cdot \delta_3 \cdot (\delta_1 + \delta_2 - n) \ge n^3/6$
  triangles in $G$, because we can pick any $w_3 \in V_3$, then any $w_1 \in N(w_3, V_1)$
  and then any $w_2 \in N(w_1) \cap N(w_3) \cap V_2$ to form a triangle.
  We will also need the following fact:
  \begin{equation}\label{eq:c3fact2}
  \text{$\forall u_1,u'_1 \in V_1$ there are at least $6 \varepsilon^2 n^2$ edges $u_2u_3$ 
  s.t.\ $u_1u_2u_3$ and $u'_1u_2u_3$ are triangles}.
  \end{equation}
  To see \eqref{eq:c3fact2}, note that there are
  at least $2 \delta_3 - n \ge 2 \varepsilon n$ ways to pick a vertex $u_3 \in V_3$ 
  adjacent to both $u_1$ and $u'_1$ and that then there
  are at least 
  $2\delta_1 + \delta_2 - 2n \ge 3 (\delta_1 + \delta_2)/2 - 2n \ge 3 \varepsilon n$
  ways to select a vertex $u_2 \in V_2$ that is adjacent to $u_1$, $u'_1$, and $u_3$.
  
  The fact that there are at least $n^3/6$ triangles and \eqref{eq:c3fact2}
  immediately implies that, for every $v,v' \in V_1$,
  the number of  $(v,v',5)$-linking sequences is at least
  $\frac{\varepsilon^3 n^5}{100}$,
  because the sequence $(u_2,u_3,w_1,w_2,w_3)$
  is a $(v,v',5)$-linking sequence whenever $vu_2u_3$ and $v'u_2u_3$ are both
  triangles and $w_1w_2w_3$ is a triangle disjoint from $\{v,v',u_2,u_3\}$. 
  
  So we are left to consider the case when $v,v' \in V_i$ for $i \in \{2, 3\}$.
  Let $j \in \{2, 3\} \setminus \{i \}$ so we have that $\{i,j\} = \{2, 3\}$.  
  We can pick $u_j \in N(v) \cap N(v') \cap V_j$ in at least 
  $2\delta_2 - n \ge 2 \varepsilon n$ ways.
  Then we can pick $u_1 \in N(v) \cap N(u_j) \cap V_1$ in at least $\delta_1 + \delta_3 - n \ge n/6$ ways.
  Similarly, we can now pick $u'_1 \in N(v') \cap N(u_j) \cap V_1$ distinct from $u_1$ 
  in at least $\delta_1 + \delta_3 - n - 1 \ge n/6$ ways.
  Observe that $vu_ju_1$ and $v'u_ju'_1$ are both triangles.
  By \eqref{eq:c3fact2}, there are at least $\frac12 \cdot 6 \varepsilon^2 n^2$ ways to
  now pick $u_2$ and $u_3$ such that $u_1u_2u_3$ and $u'_1u_2u_3$ are both triangles and
  such that $u_2$ and $u_3$ are disjoint from $\{v, v', u_j\}$.
  All together there are at least 
  \begin{equation*}
     2 \varepsilon \cdot \frac16 \cdot \frac16 \cdot 3 \varepsilon^2 \cdot n^5
     \ge \frac{\varepsilon^3 n^5}{100}
  \end{equation*}
  ways to make these selection. To complete the proof, 
  we observe that every such selection
  $(u_j, u_1, u_2, u_3, u'_1)$ is a $(v,v',5)$-linking sequence, 
  because $vu_ju_1$ and $u'_1u_2u_3$ are both triangles and
  $v'u_ju'_1$ and $u_1u_2u_3$ are both triangles.
\end{proof}  

\begin{lemma}\label{lem:weakerstrongtostrong}
  Let $k \ge 3$. 
  If Conjecture~\ref{conj:weakerstrong} holds for $k$,  
  then Conjecture~\ref{conj:strong} holds for $k$.
\end{lemma}
\begin{proof}
  We can assume $\sigma$ is small enough and $n$ is large enough
  so that the following holds:
  \begin{itemize}
    \item 
      $\sigma^{1/2} < \varepsilon$ and 
      for every $n' \ge (1-\sigma^{1/2})n$
      and for every $3 \le \ell \le k$, we can apply Conjecture~\ref{conj:weakerstrong} 
      with $\ell$, $\sigma$, $n'$ and $\varepsilon - \sigma^{1/2}$ playing the roles
      of $k$, $\sigma$, $n$ and $\varepsilon$, respectively;
    \item for every $4 \le \ell \le k$, 
      we can apply Lemma~\ref{lem:absorb} with
      $\ell$, $\ell-1$, $\varepsilon^3/2^\ell$, and $\sigma^{1/2}$
      playing the roles of $k$, $t$, $\eta$, and $\sigma$, respectively; and
    \item
      we can apply Lemma~\ref{lem:absorb} with
      $3$, $5$, $\varepsilon^3/100$, and $\sigma^{1/2}$
      playing the roles of $k$, $t$, $\eta$, and $\sigma$, respectively.
  \end{itemize}
  Let $G$ and $\delta_1, \dotsc, \delta_k$
  be as in the statement of Conjecture~\ref{conj:strong}.
  Let $I := \{i \in [k] : \delta_i < (1 + \varepsilon) \frac n2\}$, 
  let $\ell = k - |I|$, and
  let $i_1 < \dotsm < i_{|I|}$ be an ordering of the elements of $I$.
  In the manner described after the statement
  of Theorem~\ref{thm:strongk3},
  iteratively, for $j$ from $1$ to $|I|$,
  we can match every $v \in V_{i_j}$ to a unique $f_v \in V_{i_j+1}$
  and then collapse the edge $vf_v$ into $f_v$.
  Let $G'$ be the resulting graph, so
  $G'$ will be a subgraph of the $n$-blow-up of $C_\ell$
  such that a transversal $C_\ell$ factor of $G'$ corresponds
  to a transversal $C_k$ factor of $G$.
  For convenience, we relabel the parts of $G'$ as 
  $V'_1,\dotsc,V'_\ell$ so that, 
  for $i \in [\ell]$, we have $G'[V'_i, V'_{i+1}] \ge \delta'_i$.
  Note that $\delta'_i \ge (1 + \varepsilon) \frac n2$ for $i \in [\ell]$ and 
  \begin{equation}\label{eq:reduction}
    \sum_{i=1}^{\ell} \delta'_i = 
    \sum_{i=1}^{k} \delta_i  - \sum_{j=1}^{k-\ell}\delta_{i_j} 
    > k\left(1 + \frac1k + \varepsilon\right)\frac n2 - (k - \ell)(1 + \varepsilon) \frac n2 
    = \ell\left(1 + \frac 1\ell + \varepsilon\right) \frac n2.
  \end{equation}
  Clearly \eqref{eq:reduction} 
  implies $\ell \ge 2$ and that we can assume $\ell \ge 3$ by Remark~\ref{remark:C2}.  
  If $\ell=3$, then Proposition~\ref{prop:c3linked} implies that
  $G$ is $(\varepsilon^3/100,5)$-linked and if 
  $\ell \ge 4$, Proposition~\ref{prop:cklinked} implies
  that $G$ is $(\varepsilon^3/2^\ell, \ell - 1)$-linked.
  So by the selection of $\sigma$ and $n$, 
  we can apply Lemma~\ref{lem:absorb} with $G'$ and $\sigma^{1/2}$
  playing the roles of $G$ and $\sigma$ to find a 
  set $A \subseteq V(G')$ with 
  $z = |V'_i \cap A| \le \sigma^{1/2}$ for $i \in [\ell]$
  guaranteed by Lemma~\ref{lem:absorb}. 
  Conjecture~\ref{conj:weakerstrong} then implies that $G' - A$
  has a transversal $C_\ell$-tiling of size at least $n - z - \sigma n$
  which implies that $G'$ has a transversal $C_\ell$-factor.
  This in turn implies that $G$ has a transversal $C_k$-factor.
\end{proof}


\section{Proof of Theorem~\ref{thm:asymp}}

Informally the proof of Theorem~\ref{thm:asymp} proceeds as follows:
Given a spanning subgraph of the $n$-blow-up of $C_k$ with parts 
$V_1, \dotsc, V_k$
that satisfies the degree condition \eqref{eq:asymp_degree_condition}, 
we independently select, 
for every $i \in k$ and for large $T := T(k, \varepsilon)$,
a partition of almost all of $V_i$ into $T+1$ parts
$U_{i,1}, W_{i,1}, W_{i,2}, \dotsc, W_{i,T}$ each of size $mk$.
The Chernoff and union bounds imply that, if $n$ is sufficiently large,
there exists an outcome where, for every $i \in [k]$ and  
every $v \in V_{i-1} \cup V_{i+1}$, the vertex $v$
has at least $(1 + 1/k + \varepsilon/2)mk/2$ neighbors in each of
the $T+1$ parts of $V_i$.
Therefore, for  $t$ from $1$ to $T$, 
we can iteratively apply the following lemma  (Lemma~\ref{lem:asymp})
to find a transversal $C_k$-tiling $\cT_t$ of size $mk$ contained in $\bigcup_{i \in k} \left(U_{i,t} \cup W_{i,t}\right)$
so that, if, for $i \in [k]$,  we let $U_{i,{t+1}}$ be the vertices in $U_{i,t} \cup W_{i,t}$ 
uncovered by $\cT_t$, we can continue with the next iteration.
In this way, we can cover almost all of the vertices, so with absorbing 
(i.e.\ Proposition~\ref{prop:cklinked} and Lemma~\ref{lem:absorb}) we can 
find a transversal $C_k$-factor.

\begin{lemma}\label{lem:asymp}
  For $\varepsilon>0$ and integer $k \ge 3$, there exists $m_0 := m_0(k,\varepsilon)$ such
  for every $m \ge m_0$ the following holds for $n \ge 2mk$.
  Suppose that $G$ is a subgraph of the $n$-blow-up of $C_k$ with parts $V_1, \dotsc, V_k$,
  and that, for every $i \in [k]$, there exist disjoint 
  $U_i, W_i \subseteq V_i$
  where $|U_i| = |W_i| = mk$ and the following conditions hold
  for every $v \in V_{i-1} \cup V_{i+1}$: 
  \begin{enumerate}[label=(C\arabic*)]
    \item\label{C1} $d(v, U_i) \ge \left(1 + \sigma\right)mk/2$, and
    \item\label{C2} $d(v, W_i) \ge \left(1 + 1/k + \sigma\right)mk/2$.
  \end{enumerate}
  Then $G$ contains a transversal $C_k$-tiling $\cT$ of size $mk$
  contained in $\bigcup_{i \in [k]} U_i \cup W_i$
  such that for every $i \in [k]$ and every $v \in V_{i-1} \cup V_{i+1}$ with
  $U'_i := \left(U_i \cup W_i\right) \setminus V(\cT)$
  we have
  \begin{equation*}
    d(v, U'_i) \ge \left(1 + \sigma\right)mk/2.
  \end{equation*}
\end{lemma}
\begin{proof}
  For every $i \in [k]$, independently and uniformly at random select a partition of $U_i$
  into parts $U_{i,1}, \dotsc, U_{i,k}$ each of size $m$.
  For every $i,j \in [k]$ and every $v \in V_{i-1} \cup V_{i+1}$,
  the random variable $d(v, U_{i,j})$ is
  hypergeometrically distributed with expected value 
  $d(v, U_i) \frac{|U_{i,j}|}{|U_i|} = \frac{d(v, U_i)}{k}$.
  Therefore, by \ref{C1} and the Chernoff and union bounds, there exists an outcome such that
  for every $i,j \in [k]$ and every $v \in V_{i-1} \cup V_{i+1}$ we have
  \begin{equation}\label{eq:after_random_selection}
    d(v, U_{i,j}) \ge m/2.
  \end{equation}
  This implies that 
  for every $i,j\in [k]$, the bipartite graph 
  $G[U_{i,j}, U_{i+1,j}]$ is balanced with parts of size $m$ and
  minimum degree at least $m/2$, so, by Hall's Theorem,
  it contains a perfect matching $M_{i,j}$.
  For $j \in [k]$, let $H_j$ be the graph with vertex set $U_{1,j} \cup \dotsm \cup U_{k,j}$
  such that 
  \begin{equation*}
    E(H_j):=\bigcup_{i=1}^{k} M_{i,j} \setminus \left(M_{j-1,j} \cup M_{j,j}\right).
  \end{equation*}
  Note that $H_j$ consists of a collection 
  $\cP_{j}$ of $m$ vertex disjoint paths each on $k-1$ vertices such that
  \begin{itemize}
    \item $V(\cP_j) = V(H_j) \setminus U_{j,j}$;
    \item every $P \in \cP_j$ has exactly one vertex in each of the sets $U_{1,j}, \dotsc, U_{k,j}$ except
      $U_{j,j}$; and
    \item every $P \in \cP_j$ has one end-vertex in $U_{j-1,j}$ and the other end-vertex in $U_{j+1,j}$.
  \end{itemize}
  By \ref{C2},
  the number of common neighbors in $W_j$
  of the endpoints of every path in $\cP_j$ is at least 
  \begin{equation*}
    2\left(1 + 1/k + \sigma\right)mk/2 - mk > m = |P_j|.
  \end{equation*}
  Therefore, we can greedily select such a common neighbor for every path in $\cP_j$ to form 
  $\cT_j$ a transversal $C_k$-tiling of size $m$. 
  The union $\cT := \cT_1 \cup \dotsm \cup \cT_k$ is a transversal $C_k$-tiling of $G$ of size $mk$.
  For every $i \in [k]$, let
  $U'_i := \left(U_i \cup W_i\right) \setminus V(\cT) = U_{i,i} \cup \left(W_i \setminus V(\cT)\right)$, 
  so
  \begin{equation*}
    |U'_i| = |U_{i,i}| + |W_i| - m = mk. 
  \end{equation*}
  With \ref{C2} and \eqref{eq:after_random_selection},
  for every $v \in V_{i-1} \cup V_{i+1}$, we have that 
  \begin{equation*}
    d(v, U'_i) \ge m/2 + \left(1 + 1/k + \sigma\right)mk/2 - m  = \left(1 + \sigma\right)mk/2. \qedhere
  \end{equation*}
\end{proof}

\begin{proof}[Proof of Theorem~\ref{thm:asymp}]
  Define
  $\eta := k^{-3} \cdot 2^{-k}$ and $\sigma := \min\{\varepsilon/4, 0.1 \eta^{k+1}/(k^4 + 1)\}$.
  By Proposition~\ref{prop:cklinked}, $G$ is $(\eta, k-1)$-linked.
  Let $A$ be the set guaranteed by Lemma~\ref{lem:absorb},
  so there exists $z \le \sigma n$ such that $|A \cap V_i| = z$ for every $i \in [k]$.
  Let $m := \floor{\frac{\sigma^2 n}{2k}}$
  and let $T := \floor{\frac{n - z}{mk}}-1$ and
  note that $(T+1)mk \le n - z \le (T+2)mk$ and that $T$ is bounded above by
  a constant that depends only on $k$ and $\varepsilon$.
  For every $i \in [k]$, let $V'_i \subseteq V_i \setminus A$ where 
  $|V'_i| = (T+1)mk$.
  We will construct $T$ disjoint transversal $C_k$-tilings each of of 
  size $mk$ that each avoid $A$. 
  Because,
  \begin{equation*}
    mkT = (T+2)mk - 2mk \ge n - z - \sigma^2 n 
  \end{equation*}
  this will imply the theorem by the properties of $A$ from Lemma~\ref{lem:absorb}.

  Note that, by \eqref{eq:asymp_degree_condition}, for every $i \in [k]$ and $v \in V_{i-1} \cup V_{i+1}$, we have 
  \begin{equation}\label{eq:nprime_degree}
    d(v, V_i') \ge \delta^*(G) - (n - |V'_i|) \ge 
    \left(1 + 1/k + 2\sigma\right)n/2 \ge \left(1 + 1/k + 2\sigma\right)|V_i'|/2.
  \end{equation}
  For every $i \in [k]$, independently and uniformly at random
  select a partition of $V'_i$ into $T+1$ parts $W_{i,0}, \dotsc, W_{i,T}$ each of size $mk$.
  For every $i \in [k]$, every $0 \le t \le T$, and every $v \in V_{i-1} \cup V_{i+1}$,
  the random variable
  $d(v, W_{i, t})$
  is hypergeometrically distributed with expected value 
  $d(v, V'_i) \frac{|W_{i,t}|}{|V'_i|}$.
  Therefore, by \eqref{eq:nprime_degree} and the Chernoff and union bounds, 
  there exists an outcome such that
  for every $i \in [k]$, $0 \le t \le T$, and $v \in V_{i-1} \cup V_{i+1}$ we have
  \begin{equation}\label{eq:after_random_selection_vW} 
    d(v, W_{i, t}) \ge \left(1 + 1/k + \sigma\right)mk/2.
  \end{equation}

  We will now show by induction on $t$ from $1$ to $T+1$ that there exist $t-1$ disjoint transversal $C_k$-tilings 
  $\cT_1, \dotsc, \cT_{t-1}$ each of size $mk$ that are contained in $\bigcup_{i \in k} \bigcup_{s = 0}^{t-1} W_{i,s}$,
  and that, for every $i \in [k]$, 
  if we let $U_{i, t} := \left(\bigcup_{s=0}^{t-1} W_{i,s}\right) \setminus \left(\bigcup_{s = 1}^{t-1} V(\cT_s)\right)$,
  then the following holds:
  \begin{equation}\label{eq:Ut_condition}  
    d(v, U_{i, t}) \ge (1 + \sigma)mk/2 \qquad \text{ for every $v \in V_{i-1} \cup V_{i+1}$.}
  \end{equation}
  This will prove the theorem.
  
  For the base case, note that when $t = 1$ we have that $U_{i,t} = U_{i,1} =  W_{i,0}$
  for every $i \in [k]$ 
  so \eqref{eq:Ut_condition} holds by \eqref{eq:after_random_selection_vW}.
  Now assume the induction hypothesis holds for some $1 \le t \le T$.
  With \eqref{eq:after_random_selection_vW} and \eqref{eq:Ut_condition}
  we can apply Lemma~\ref{lem:asymp} 
  to find a tiling $\cT_t$ of size $mk$ contained in 
  $\bigcup_{i \in [k]} U_{i, t} \cup W_{i, t}$ such that,
  for every $i \in [k]$, we have that 
  \begin{equation*}
    U_{i,t+1} = 
    \left(\bigcup_{s=0}^{t} W_{i,s}\right) \setminus \left(\bigcup_{s = 1}^{t} V(\cT_s)\right) =
    \left(U_{i, t} \cup W_{i, t}\right)\setminus V(\cT_t)
  \end{equation*} 
  satisfies \eqref{eq:Ut_condition} with $t$ set to $t+1$.  Therefore, the induction hypothesis holds for $t+1$.
\end{proof}

\section{Proof of Theorem~\ref{thm:strongk3}}

Because Theorem~\ref{thm:almoststrongk3} works for every $n$ and the degree
condition is weaker than Theorem~\ref{thm:strongk3}, 
it might have independent interest.
Note that Theorem~\ref{thm:almoststrongk3} is stronger than 
the $k=3$ case of Conjecture~\ref{conj:weakerstrong}, so
Lemma~\ref{lem:weakerstrongtostrong} and Theorem~\ref{thm:almoststrongk3}
together imply Theorem~\ref{thm:strongk3}.

\begin{theorem}\label{thm:almoststrongk3}
  The following holds for every $n \in \mathbb{N}$ and
  every subgraph $G$ of the $n$-blow-up of $C_3$ with parts $V_1, V_2, V_3$.
  If there exist $\delta_1, \delta_2,\delta_3\ge n/2$ such that
  $\delta_1 + \delta_2 + \delta_3 \ge 2n$ and
  \begin{equation*}
    \delta(G[V_i, V_{i+1}]) \ge \delta_i \qquad\text{for every $i \in [3]$},
  \end{equation*}
  then $G$ has a transversal $C_3$-tiling of size at least $n-1$.
\end{theorem}
\begin{proof}
  For brevity, in this proof we call a transversal $C_3$-tiling
  a \textit{tiling}.
  Let the size of a maximum tiling of $G$ be $m$ and let us assume for a 
  contradiction that $m \le n-2$.
  Call a pair of edges $e$ and $f$ \textit{dissimilar} if 
  $e \in E(G[V_i, V_{i+1}])$ and $f \in E(G[V_j, V_{j+1}])$ for
  distinct $i,j \in [3]$.
  Call a set $F \subseteq E(G)$ a \textit{dissimilar matching} if 
  the edges in $F$ are disjoint and the edges in $F$ are pairwise
  dissimilar.
  For every maximum tiling $\cT$, let $h(\cT)$ be the maximum size of
  a dissimilar matching $F \subseteq E(G[U(\cT)])$ such that
  every edge in $F$ is uncovered by $\cT$.
  Recall that an edge $e$ is uncovered by $\cT$ if both endpoints of
  $e$ are disjoint from $V(\cT)$.
  Let $\{\alpha, \beta, \gamma\} = \{\delta_1/n, \delta_2/n, \delta_3/n\}$ 
  and $\{A, B, C\} = \{V_1, V_2, V_3\}$ 
  be labellings such that $\alpha \le \beta \le \gamma$ and 
  \begin{equation}\label{eq:abccond}
    \delta(G[B,C]) \ge \alpha n, \quad
    \delta(G[A,C]) \ge \beta n, \quad\text{and}\quad
    \delta(G[A,B]) \ge \gamma n.
  \end{equation}

  \begin{claim}\label{clm:dissimilar}
    Let $\cT$ be a maximum tiling.
    If $e \in E(G)$ is uncovered
    by $\cT$, then $d(e, U(\cT)) = 0$.
    Furthermore, if $e$ and $f$ are disjoint dissimilar edges that are 
    uncovered by $\cT$, then
    $d(e, T) + d(f, T) \le 1$ for every $T \in \cT$.
  \end{claim}
  \begin{proof}
    If $e$ is an edge uncovered by $\cT$ and $x \in U(\cT)$ is such that
    $d(e, \{x\}) = 1$, then $ex$ is triangle, and adding $ex$ to $\cT$
    creates a tiling of size $m+1$, a contradiction.
    Similarly, if $e$ and $f$ are disjoint and dissimilar edges
    that are uncovered by $\cT$ and $d(e, T) + d(f, T) \ge 2$
    for some $T \in \cT$, then,
    because $e$ and $f$ are dissimilar, there exist distinct
    $x,y \in T$ such that $ex$ and $fy$ are both triangles, so
    if we replace $T$ with $ex$ and $fy$ in $\cT$, then we have
    a tiling of size $m+1$, a contradiction.
  \end{proof}

  \begin{claim}\label{clm:hle2}
    Let $\cT$ be a maximum tiling 
    and let $F$ be a dissimilar matching with $|F|=3$.
    Then either there exists $e \in F$ such that $d(e, U(\cT)) \ge 1$
    or there exists $T \in \cT$ such that
    $\sum_{e \in F} d(e, T) \ge 2$.
    Consequently, $h(\cT) \le 2$ for every maximum tiling $\cT$.
  \end{claim}
  \begin{proof}
    Let $\{a_1, \dotsc, a_n\}$, $\{b_1, \dotsc, b_n\}$, and $\{c_1, \dotsc, c_n\}$
    be orderings of $A$, $B$, and $C$, respectively, such that 
    $a_ib_ic_i \in \cT$ for every $i \in [m]$ (so, when $m + 1 \le i \le n$, 
    $a_ib_ic_i$ is not a triangle).
    By \eqref{eq:abccond}, we have     
    \begin{equation*}
      \sum_{i=1}^{n} \sum_{e \in F} d(e, a_ib_ic_i)=
      \sum_{e \in F} d(e, V(G)) \ge \\
      (\alpha + \beta - 1)n + (\alpha + \gamma - 1)n + (\beta + \gamma - 1)n \ge n.
    \end{equation*}
    Therefore, if $0=\sum_{e \in F} d(e, U(\cT)) = 
    \sum_{i=m+1}^{n}\sum_{e \in F} d(e, a_ib_ic_i)$, then 
    there exists 
    $1 \le i \le m$ such that $\sum_{e \in \cF} d(e, a_ib_ic_i) = 2$.
    This proves the first statement. 

    To see the second statement, assume that for a contradiction
    that there exists a maximum tiling $\cT$ such that $h(\cT) = 3$.
    This means that there exists a dissimilar matching  
    $F$ such that $|F| = 3$ and such that every edge in $F$ is uncovered
    by $\cT$.  By the first part of the statement, there either exists $e \in F$ such  
    $d(e, U(\cT)) \ge 1$, or there exists two edges $e,f \in \cT$ and $T \in \cT$
    such that $d(e, T) + d(f, T) \ge 1$.
    Because every edge in $F$ is uncovered by $\cT$, 
    this contradicts Claim~\ref{clm:dissimilar}.
  \end{proof}

  \begin{claim}\label{clm:edgeswap}
    There exists a maximum tiling $\cT$ such that $h(\cT) = 2$,
    and for every maximum tiling $\cT$ there does not exist $e\in E(G[B,C])$ 
    which is uncovered by $\cT$.
  \end{claim}
  \begin{proof}
    Suppose for a contradiction that the statement is false and
    assume $\cT$ and a dissimilar matching $F$ in $G[U(\cT)]$ have both
    been selected so that 
    \begin{enumerate}[label=(\Alph*)]
      \item\label{A} there exists $e \in F$ such that $e \in G[B,C]$ if possible, and,
      \item\label{B} subject \ref{A}, $|F|$ is as large as possible.
    \end{enumerate}
    Note that Claim~\ref{clm:hle2}, implies that $|F| \le h(\cT) \le 2$, 
    so if there exists $e \in F$ such that $e \in E(G[B,C])$, then 
    $F$ has at most one edge that is contained in $E(G[A,B]) \cup E(G[A,C])$.
    If there is no $e \in F$ that is in $E(G[B,C])$, then 
    by the selection of $\cT$ and $F$ (c.f.\ \ref{A}),
    for every maximum tiling $\cT$ there does not exist $e \in E(G[B,C])$,
    so our contrary assumption implies $|F| \le h(\cT) \le 1$.
    Therefore, in all cases, 
    $F$ has at most one edge that is contained in $E(G[A,B]) \cup E(G[A,C])$.
    Let $\{X, Y\} = \{B, C\}$ be a labelling such that
    $F$ does not contain an edge in $G[A,X]$.

    Let $W \subseteq U$ be the set of vertices that are incident to an edge
    in $F$.
    The fact that $|\cT| \le n -2$,  implies that there exist nonadjacent vertices 
    $a \in A \setminus W$ and $x \in X \setminus W$ that are uncovered by $\cT$.
    Let $\{a_1, \dotsc, a_n\}$, $\{x_1, \dotsc, x_n\}$ and 
    $\{y_1, \dotsc, y_n\}$ be orderings of $A$, $X$, and $Y$, respectively
    such that $a_n = a$, $x_n = x$ and 
    $a_ix_iy_i \in \cT$ for every $i \in [m]$.
    We can assume that the orderings are such that $W$
    is contained in the set $\{a_{n-1}, x_{n-1}, y_{n-1}, y_n\}$ with
    $x_{n-1}y_{n-1} \in F$ if $F \cap E(G[X,Y]) = F \cap E(G[B,C]) \neq \emptyset$.

    Since $a$ and $x$ are nonadjacent,
    $d(x, a_n) + d(a, x_n) + d(a, y_n) = d(x, a) + d(a, x) + d(a, y_n) \le 1$,
    and, by \eqref{eq:abccond} and the fact that $\alpha \le \beta \le \gamma$, we have
    \begin{equation*}
      \sum_{i=1}^{n} d(x, a_i) + d(a, x_i) + d(a, y_i) = 
      d(x, A) + d(a, X) + d(a, Y) \ge \beta n + \beta n + \gamma n \ge 2n,
    \end{equation*}
    so there must exist $i \in [n-1]$ such that 
    $d(x, a_i) + d(a, x_i) + d(a, y_i) = 3$.
    Note that $i \neq n-1$, because if $ax_{n-1}$ is an edge, then the
    fact that $F \cap E(G[A,X]) = \emptyset$ and the maximality of $F$
    imply that 
    $x_{n-1}y_{n-1} \in F$, but, because $\cT$ is a maximum tiling,
    $ax_{n-1}y_{n-1}$ is not a triangle.
    Since $F \cap E(G[A,X]) = \emptyset$, the maximality of $F$ also
    implies that that 
    $d(x, a_j) = 0$ for every $m+1 \le j \le n-2$, so 
    it must be that $i \in [m]$, i.e., that 
    $T = a_ix_iy_i$ is a triangle in $\cT$. 
    Therefore, we can swap $T$ for the triangle $ax_iy_i$ in $\cT$
    to form the maximum tiling $\cT'$. 
    Because the edge $a_ix$ is uncovered by $\cT'$ and $W \subseteq U(\cT')$, 
    we have a contradiction to the selection of $\cT$ and $F$ (c.f.\ \ref{B}).
  \end{proof}

    Claim~\ref{clm:edgeswap} implies 
    that there exists a maximum tiling
    $\cT$ such that $h(\cT) = 2$.
    By Claim~\ref{clm:edgeswap}, 
    we can assume that $\cT$ leaves no
    edge in $G[B,C]$ uncovered by $\cT$.
    This means that there are
    disjoint edges $ab \in G[A, B]$ and
    $a'c \in G[A,C]$ with $a, a' \in A$ that are uncovered by $\cT$.
    Since $|\cT| = m \le n - 2$, 
    there also exists $b' \in B \setminus \{b\}$ and 
    $c' \in C \setminus \{c\}$ that are uncovered by $\cT$.
    Furthermore, the fact that no edge in 
    $G[B,C]$ is uncovered implies that 
    \begin{equation*}
      \sum_{T \in \cT} d(b', T \cap C) + d(c', T \cap B) =
      d(b', C) + d(c', B) \ge 2 \delta(G[B,C]) \ge n > m = |\cT|,
    \end{equation*}
    so there exists $T \in \cT$ such that 
    $d(b', T \cap C) + d(c', T \cap B) = 2$.
    Let $e$ be the edge incident to $c'$ and $T \cap B$
    and let $e'$ be the incident to $b'$ and $T \cap C$.
    If we define $F := \{ab, a'c, e\}$, then $F$ is a dissimilar matching,
    so Claim \ref{clm:hle2} implies that we are in one of the following two 
    cases.

    \noindent \textbf{Case 1}:  
    \textit{There exists $f \in F$ such that $d(f, U(\cT)) \ge 1$.}
    By Claim~\ref{clm:dissimilar}, the fact that $ab$ and $a'c$ are
    uncovered by $\cT$ implies that $f = e$.
    So, there exists a triangle $T'$ that contains
    $e$ and a vertex in $U \cap A$. 
    This means that we can create a maximum tiling by replacing $T$ with $T'$
    in $\cT$ that leaves the edge $e' \in G[B,C]$ uncovered, 
    contradicting Claim~\ref{clm:edgeswap}.

    \noindent \textbf{Case 2}:
    \textit{There exists $T' \in \cT$ such that $\sum_{f \in F} d(f,T') \ge 2$.}
    By Claim~\ref{clm:dissimilar}, $d(ab, T') + d(a'c, T') \le 1$, so there
    is $f \in \{ab, a'c\}$ such that $d(e, T') + d(f, T') \ge 2$.
    This means that there exist two triangles, say $T''$, $T'''$, 
    in the graph induced by the vertices incident to $e$, $f$ and $T'$.
    Therefore, we can create a new tiling, say $\cT'$, by removing $T$ and $T'$
    from $\cT$ and replacing them with  $T''$ and $T'''$.
    Since $\cT$ is maximum tiling, we have that $T \neq T'$ and that
    $\cT'$ is a maximum tiling.
    Because $T \neq T'$, the edge $e' \in G[B,C]$ is uncovered by $\cT'$
    which contradicts Claim~\ref{clm:edgeswap}.
%
\end{proof}


\section{Proof of Theorem \ref{thm:trianglecover}}\label{sec:triangle_covers}

The following example proves the second part of the theorem
\begin{example}
For the case $\gamma=2/3$, it can be checked that Example \ref{exp:conj} for $k=3$ satisfies the second claim of Theorem \ref{thm:trianglecover}.


For the case $\gamma \in (\frac34,\frac79]$, we assume $n$ satisfies the following: $\gamma\ge \frac34+\frac1n$ and $(1-\beta)n/2$ is an integer. 
Clearly, since $\beta$ is rational and $\gamma>4/3$, there are infinitely many choices for such $n$.
Let us fix $\varepsilon\in (0,\frac1n]$ such that $(1-\gamma+\varepsilon)n$ is an integer.

Take sets $A=A_0\cup A_1\cup A_2 \cup A_3$, $B=B_0\cup B_1\cup B_2 \cup B_3$ and $C=C_0\cup C_1\cup C_2 \cup C_3$, such that: 
\begin{itemize}
  \item $|B_i|=(1-\gamma+\varepsilon)n$, $|A_i|=|C_i|=(1-\beta)n/2$, for $i\in[3]$;  

  \item $|B_0|=n-3(1-\gamma+\varepsilon)n=(3\gamma -2)n-3\varepsilon n$; and

  \item $|A_0|=|C_0|=n-3(1-\beta)n/2=(3\beta-1)n/2$.
\end{itemize}

Let $G$ be the $3$-partite graph with parts $A,B$ and $C$, where $E(G)$ consists of the union of the 
  edges in the following graphs:
  \begin{itemize}
    \item the complete bipartite graphs with parts $A_0, B\cup C$ and $B_0, A\cup C$ and $C_0, A\cup B$.
     \item the complete bipartite graphs with parts $A_1,B_2\cup B_3$ and $A_2,B_1\cup B_3$ and $A_3,B_1\cup B_2$.
     \item the complete bipartite graphs with parts $B_i,C_i$ and $A_i,C_i$ for each $i\in[3]$.
  \end{itemize}
Since $\gamma \le \frac 79$ and $\varepsilon > 0$, we have $(1-\gamma+\varepsilon)n > (\gamma - 1/3)n/2 =(1-\beta)n/2$. So, 
\begin{itemize}
  \item $\delta(G[A,B])=n-(1-\gamma+\varepsilon)n= \gamma n-\varepsilon n$,
  \item $\delta(G[B,C])= n-2(1-\gamma+\varepsilon)n= (2\gamma-1) n-2\varepsilon n$, and 
  \item $\delta(G[A,C])= n-2(1-\beta)n/2=\beta n$. 
\end{itemize}
Recall that $\varepsilon\le \frac1n$ and $\gamma\ge \frac34+\frac1n$, so  $\delta(G[A,B]) \ge \gamma n -1$ and
$\delta(G[B,C]) \ge (2\gamma - 1)n - 2 \ge n/2$. 

%
Note that $A_0 \cup B_0 \cup C_0$ is a triangle cover and
$$|A_0|+|B_0|+|C_0|=(3\beta-1)n/2+
(3\gamma -2)n-3\varepsilon n+(3\beta-1)n/2=(1-3\varepsilon)n<n.$$ 
\end{example}

We now proceed with the proof of the first part of the Theorem~\ref{thm:trianglecover}.
We will use the following definition throughout the proof.
\begin{definition}
  For $U, W, U' \subseteq V(G)$, let 
  $P_{3}(U, W, U')$ be the set of paths on $3$ vertices in which the middle
  vertex is in $W$, one endpoint is in $U$ and the other endpoint is in $U'$.
  When a set $\{u\}$ is a singleton, 
  we sometimes replace $\{u\}$ with $u$ in this notation.
\end{definition}

Let $\{\alpha, \beta, \gamma\} = \{\delta_1/n, \delta_2/n, \delta_3/n\}$ 
and $\{A, B, C\} = \{V_1, V_2, V_3\}$ 
be labellings such that $\alpha \le \beta \le \gamma$ and 
\begin{equation*}
\delta(G[B,C]) \ge \alpha n, \quad
\delta(G[A,C]) \ge \beta n, \quad\text{and}\quad
\delta(G[A,B]) \ge \gamma n.
\end{equation*}
We can assume $\gamma + \beta = 4/3$ and $\alpha = 1/2$. Therefore, 
  \begin{equation}\label{eq:betagammaa}
    1/2 \le \beta \le 2/3 \qquad\text{and}\qquad 2/3 \le \gamma \le 5/6.
  \end{equation}
  Let $U$ be a triangle cover and let $x = |A \cap U|/n$,
  $y = |B \cap U|/n$ and $z = |C \cap U|/n$.
  For a contradiction, assume that 
  \begin{equation}\label{eq:contrary_assumption}
    x + y + z < 1.
  \end{equation}
  Let $A' = A \setminus U$, $B' = B \setminus U$, and $C' = C \setminus U$.

  \begin{claim}\label{clm:1}
    $x \ge 1/3$, $y \ge \gamma - 1/2$, and $z \ge \beta - 1/2$.
  \end{claim}
  \begin{proof}
    Since $y + z \le x + y + z < 1$, one of $y$ or $z$ is less than $1/2$, so there
    exists an edge $bc \in G[B',C']$. Because $G[A', B', C']$ is triangle-free,
    \begin{equation*}
      0 = |N(b, A') \cap N(c, A')| \ge  
      d(b, A) + d(c, A) - |A| - xn \ge \gamma n + \beta n -  n  - x n = n/3 - x n,  
    \end{equation*}
    so $x \ge 1/3$.  By considering an edge in $G[A',B']$ and 
    an edge in $G[A',C']$ the same argument yields
    $z \ge \beta - 1/2$ and $y \ge \gamma - 1/2$, respectively.
  \end{proof}

  \begin{claim}\label{clm:2}
    $x < \gamma$, $y < 1/2$, and $z < 1/2$. 
  \end{claim}
  \begin{proof}
    We first show that both $x < \gamma$ and $y < \gamma$.
    To this end, note that if $x \ge \gamma$, then
    $x + y \ge \gamma + \gamma - 1/2 = 2\gamma - 1/2$.
    If $y \ge \gamma$, then, because $\gamma - 1/2 \le 5/6 - 1/2 = 1/3$, we also have 
    $x + y \ge 1/3 + \gamma \ge 2\gamma - 1/2$.
    So, in either case, we have the following contradiction 
    \begin{equation*}
      1 > x + y + z \ge 2\gamma - 1/2 + \beta - 1/2 = 1/3 + \gamma \ge 1.
    \end{equation*}
    Similarly, it is clear that $z < 1/2$, since otherwise 
    $x + y + z \ge 1/3 + (\gamma - 1/2) + 1/2 \ge 1$, a contradiction.

    Assume $y \ge 1/2$ and let $b \in B'$.
    Note that there are at most $(1-\gamma)n$ vertices
    $a \in N_{\overline{G}}(b, A')$. For every such $a$ we can find a $4$-vertex path $ab'a'b$ in $G[A',B']$. Indeed, since $y < \gamma$, there exists 
    $b' \in N(a, B')$.
    Then, because $2 \gamma + \beta \ge 2$ and $1 > x + y + z$,
    \begin{equation*}
      x < 1 - y - z \le 1/2 - z \le 1/2 - (\beta - 1/2) = 1 - \beta \le 
      2\gamma - 1.
    \end{equation*}
    Since $|N(b, A) \cap N(b', A)| \ge 2\gamma n - n>xn$, there exists
    $a' \in |N(b, A') \cap N(b', A')|$, giving us the $4$-path $ab'a'b$.
    Note that $N(a', C')$ and $N(b', C')$ are disjoint and that 
    every $c \in C'$ that is adjacent to both $a$ and $b$
    is not adjacent to $a'$ and not adjacent to $b'$.
    Therefore, $|P_3( b, C', a)|$ is at most
    \begin{equation*}
      |C' \setminus \left(N(a', C') \cup N(b', C')\right)| \le
      (1 - z)n - (\beta - z)n - (1/2 - z)n = (z - \beta + 1/2)n,
    \end{equation*}
    and $|P_3(b, C', A')|/n^2 \le (1 - \gamma)(z - \beta + 1/2)$.
    On the other hand, 
    \begin{equation*}
      |P_3(b, C', A')| \ge \sum_{c \in N(b, C')}d(c, A') 
      \ge (1/2 - z)n \cdot (\beta - x)n.
    \end{equation*}
    The claim then follows because there are no solutions to
    \begin{equation*}
      (1 - \gamma)(z - \beta + 1/2) \ge (1/2 - z)(\beta - x),
    \end{equation*}
    when $x \ge 1/3$, $y \ge 1/2$, and $z \ge \beta - 1/2$.
    (See Lemma~\ref{lem:app1} in the appendix for a proof of this fact.)
  \end{proof}

  Note that Claim~\ref{clm:2} implies that $\delta(G[A',B']) \ge 1$,
  $\delta(G[B',C']) \ge 1$, and that every vertex in $A'$ has a neighbor
  in $C'$.  (We do not yet know if every vertex in $C'$ has
  a neighbor in $A'$.)  We will use these facts in the rest of
  the argument without comment.

  In particular, 
  the fact that every $c \in C'$ has a neighbor $b \in B'$ implies that 
  \begin{equation*}
    d(c, A') \le |A' \setminus N(b, A')| \le |A \setminus N(b, A)| \le (1 - \gamma)n,
  \end{equation*}
  so 
  $|E(A', C')| = \sum_{c \in C'}d(c, A') \le |C'|(1 - \gamma)n = (1-z)(1-\gamma)n^2$.
  On the other hand, we have that 
  $|E(A', C')| = \sum_{a \in A'}d(a, C') \ge |A'|(\beta - z)n = (1-x)(\beta-z)n^2$.
  This yields the following useful inequality
  \begin{equation}\label{eq:simple}
    (1 - \gamma)(1 - z) \ge (1-x)(\beta - z).
  \end{equation}

  \begin{claim}\label{clm:A'B'path}
    For every $a \in A'$ and $b \in B'$, there is an $(a,b)$-path in $G[A', B']$
    with at most $4$ vertices.
  \end{claim}
  \begin{proof}
    Assume the contrary and let $a \in A'$ and $b \in B'$ be such that
    there is no $(a,b)$-path in $G[A', B']$ with at most $4$ vertices.
    Let $b' \in N(a, B')$ and $a' \in N(a, A')$.
    By our contrary assumption, we have 
    that $N(a, B') \cap N(a', B') = \emptyset$ so
    $y \ge |N(a, B) \cap N(a', B)|/n \ge 2\gamma - 1$,
    By the same argument, $N(b, A') \cap N(b', A') = \emptyset$
    and $x \ge 2\gamma - 1$.
    Since $\beta \le 2/3$, we have 
    $2 \gamma - 1 = 2(4/3 - \beta) - 1 = 5/3 - 2\beta \ge 1 - \beta$, so
    \begin{equation*}
      z < 1 - x - y \le 1 - 2(2 \gamma - 1) \le 
      1 - 2(1 - \beta) = 
      2 \beta - 1 \le |N(a, C) \cap N(a', C)|/n,
    \end{equation*}
    therefore there exists $c \in N(a, C') \cap N(a', C')$.
    But then $N(c, B')$ cannot intersect $N(a, B') \cup N(a', B')$, so
    \begin{equation*}
      (1/2 - y) + 2(\gamma - y) \le |N(c, B') \cup N(a, B') \cup N(a', B')|/n \le 1 - y
    \end{equation*}
    which implies that $y \ge \gamma - 1/4$, therefore $y\ge 5/12$. But \eqref{eq:simple} has no
    solutions when $y \ge 5/12$, $x \ge 1/3$ and $z \ge \beta -1/2$.
    (See Lemma~\ref{lem:app2} in the appendix for a proof of this fact.) This is a contradiction.
  \end{proof}

  \begin{claim}\label{clm:yltonethird}
    $y < 1/3$ and $x<\beta$.
  \end{claim}
  \begin{proof}
    Let $a \in A'$. 
    We first get an upper-bound on $|P_3(a, C', B')|$.
    Note that there are at most $(1 - \gamma)n$ ways to select
    $b \in B'$ that is not adjacent to $a$.  By Claim~\ref{clm:A'B'path},
    there exists $a' \in A'$ and $b' \in B'$ such that $ab'a'b$ is a path.
    Note that every vertex $c \in C'$ that is adjacent to both $a$ and $b$
    cannot be in $N(a', C') \cup N(b', C')$.
    Since $N(a', C')$ and $N(b', C)$ are disjoint,
    we have the cardinality of $P_3(a, C', b)$ is at most 
    \begin{equation*}
      |C'| - d(a', C') - d(b', C') \le 
      (1-z)n - (\beta - z)n - (1/2 - z)n = (z - \beta + 1/2)n
    \end{equation*}
    Therefore, $|P_3(a,C', B')|/n^2 \le (1 - \gamma)(z - \beta + 1/2)$.
    We also have that 
    $|P_3(a, C', B')| \ge \sum_{c \in N(a, C')}d(c, B')\ge (\beta - z)n(1/2 - y)n$,
    so
    \begin{equation}\label{eq:m1}
      (1 - \gamma)(z - \beta + 1/2) \ge 
      |P_{3}(a,C', B')|/n^2 \ge (\beta - z)(1/2 - y).
    \end{equation}
    
   By considering $b \in B'$ and estimating $P_{3}(A',C', b$), the same arguments yield that
    \begin{equation}\label{eq:m2}
      (1 - \gamma)(z - \beta + 1/2) \ge 
      |P_{3}(A',C', b)|/n^2 \ge \sum_{c \in N(b, C')}d(c,A')/n^2 \ge (\beta - x)(1/2 - z).
    \end{equation}
    But \eqref{eq:m1}, \eqref{eq:m2} and \eqref{eq:simple} cannot hold simultaneously
    when $x \ge 1/3$, $y \ge 1/3$ and $z \ge \beta - 1/2$.
    (See Lemma~\ref{lem:app3} in the appendix for a proof of this fact.)
    Therefore $y < 1/3$. 
    
    Now we will show that $x<\beta$. Indeed, if $\beta\le x$, we have 
    \begin{equation*}
        y\ge \gamma-1/2= 5/6-\beta\ge 5/6-x > y+z-1/6,
    \end{equation*} 
    so $z < 1/6$. With \eqref{eq:m1} we get 
    $(\beta-1/3)(1/6 - \beta + 1/2) \ge (\beta - 1/6)(1/2 - y)$.
    Plugging $y<1/3$ we get that $-\beta^2+(5/6)\beta-1/4>0$ which does not have a solution, a contradiction.
  \end{proof}
  Note that Claims~\ref{clm:2} and ~\ref{clm:yltonethird} together
  imply $\delta(G[A',B']), \delta(G[B',C']), \delta(G[C',A']) \ge 1$.

  \begin{claim}\label{clm:4path}
    There exists $a_1 \in A'$ and $c_1 \in C'$ such that 
    there is no $(a_1, c_1)$-path in $G[A', C']$ with at most $4$-vertices.
  \end{claim}
  \begin{proof}
    Assume the contrary and let $a_1 \in A'$.
    Then, for every $c_1 \in C'\setminus N(a',C')$, there exists
    $a_2 \in A'$ and $c_2 \in C'$ such that $a_1c_2a_2c_1$ is a path,
    so, since $G[A',B',C']$ is triangle-free,
    $|P_3(a_1, B', c_1)|$ is at most 
    \begin{equation*}
      |B' \setminus (N(a_2, B') \cup N(c_2, B')| \le
      |B'| - (d(a_2, B) - yn + d(c_2, B) - yn) \le  
      \left(y - \gamma + 1/2\right)n.
    \end{equation*}
    Since $a_1$ has at most $(1 - \beta)n$ non-neighbors in $C'$,
    we have that
    \begin{equation*}
      (1 - \beta)\left(y - \gamma + 1/2\right) 
      \ge |P_3(a_1, B', C')|/n^2 = \sum_{b \in N(a_1, B')}d(b, C')/n^2 \ge
      (\gamma - y) \left(1/2 - z\right)
    \end{equation*}
    which is impossible when 
    $x \ge 1/3$, $1/3 > y \ge \gamma - 1/2$, and $z \ge \beta - 1/2$. 
    (See the Lemma~\ref{lem:app4} in the appendix for a proof of this fact.)
  \end{proof}

  By Claim~\ref{clm:4path}, there exists $a_1 \in A'$ and $c_1 \in C'$
  such that there is no $(a_1,c_1)$-path 
  in $G[A', C']$ with at most $4$-vertices.
  Fix such vertices $a_1$ and  $c_1$.  By Claims \ref{clm:2} and \ref{clm:yltonethird},
  we can also fix 
  $c_2 \in N(a_1, C')$ and $a_2 \in N(c_1, A')$.
  Note that, by the selection of $a_1$ and $c_1$, 
  \begin{equation}\label{eq:disjoint}
    \text{$N(a_1, C') \cap N(a_2, C') = \emptyset$ and 
    $N(c_1, C') \cap N(c_2, C') = \emptyset$.}
  \end{equation}

  \begin{claim}\label{clm:largeZ}
    $z \ge \beta - 1/4$.
  \end{claim}
  \begin{proof}
    Since $|N(a_1, B) \cap N(a_2, B)|/n \ge 2 \gamma - 1 \ge 1/3 > y$,
    there exists $b \in B'$ that is adjacent to both $a_1$ and $a_2$.
    Since $G[A',B',C']$ is triangle-free \eqref{eq:disjoint} implies that 
    \begin{equation*}
      1/2 - z \le d(b, C')/n \le 
      |C' \setminus (N(a_1, C') \cup N(a_2, C'))|/n \le 
      1 - z - 2(\beta - z) = 1 - 2\beta + z,
    \end{equation*}
    so $z \ge \beta - 1/4$.
  \end{proof}

  \begin{claim}\label{clm:3disjoint}
    At least one of the following statements is true.
    \begin{itemize}
      \item
        For every $a \in A'$, we have that
        $N(a, C')$ intersects $N(a_1, C') \cup N(a_2, C')$.
      \item
        For every $c \in C'$, we have that
        $N(c, A')$ intersects $N(c_1, A') \cup N(c_2, A')$.
    \end{itemize}
  \end{claim}
  \begin{proof}
    Assume the contrary, so there exists $a_3 \in A'$ such that
    $N(a_1,C')$, $N(a_2, C')$, and $N(a_3, C')$ are pairwise disjoint and 
    that there exists $c_3 \in C'$ such that
    $N(c_1,A')$, $N(c_2, A')$, and $N(c_3, A')$ are pairwise disjoint.
    This implies that 
    \begin{equation*}
      (1-x)n = |A'| \ge d(c_1, A') + d(c_2, A') + d(c_3, A') \ge 3(\beta - x)n,
    \end{equation*}
    so $x \ge (3 \beta - 1)/2$, and, by considering the sets $N(a_1, C')$,
    $N(a_2, C')$ and $N(a_3, C')$, 
    we similarly have that $z \ge (3\beta - 1)/2$.
    This implies that
    \begin{equation*}
      y < 1 - x - z \le 1 - (3\beta - 1) = 2 - 3\left(4/3 - \gamma\right) 
      = 3\gamma - 2.
    \end{equation*}
    Note that $|N(a_1, B) \cap N(a_2, B) \cap N(a_3, B)| \ge 3 \gamma n - 2|B|
    = (3 \gamma - 2)n > y n$,
    so there exists $b \in N(a_1, B') \cap N(a_2, B') \cap N(a_3, B')$.
    Note that $N(b, C')$ must be disjoint from 
    $N(a_1, C') \cup N(a_2, C') \cup N(a_3, C')$ 
    so, since $N(a_1, C')$,$N(a_2, C')$, and $N(a_3, C')$ are pairwise disjoint,
    \begin{equation*}
      (1 - z)n = |C'| \ge d(b, C') + d(a_1, C') + d(a_2, C') + d(a_3, C')
     \ge (1/2 - z + 3(\beta - z))n,
    \end{equation*}
    so $z \ge \beta - 1/6$.
    But then $1 > x + y + z \ge 1/3 + \gamma - 1/2 + \beta - 1/6 = 1$, a contradiction.
  \end{proof}

  \begin{claim}\label{clm:ac_ipath}
    For every $a \in A'$ there exists $i \in \{1, 2\}$ such that there 
    is an $(a, c_i)$-path in $G[A',C']$ with at most $4$ vertices.
  \end{claim}
  \begin{proof}
    Since $c_1a_2$ and $c_2a_1$ are edges, we have the desired path if 
    $N(a, C')$ intersects either
    $N(a_1, C')$ or $N(a_2, C')$.  So assume otherwise, i.e.,
    assume that the sets 
    $N(a, C')$, $N(a_1, C')$, and $N(a_2, C')$ are pairwise disjoint.
    By Claim~\ref{clm:2}, there exists $c \in N(a, C')$.
    Because $N(c_1, A')$ and $N(c_2, A')$ are disjoint,
    Claim~\ref{clm:3disjoint} implies that $N(c, A')$ must intersect
    one of $N(c_1, A')$ or $N(c_2, A')$ and this gives us the
    desired path.
  \end{proof}

  For $i \in \{1,2\}$, 
  let $A_i = N(c_i, A')$, 
  let $B_i \subseteq N(c_i, B')$ such that $|B_i| = \ceiling{(1/2 - y)n}$,
  and let $B_0 = B' \setminus (B_1 \cup B_2)$.
  Define $\zeta = |B_1|/n = |B_2|/n$, so
  $|B_0| = (1 - y - 2\zeta)n$.

  \begin{claim}\label{clm:path}
    Every $a \in A_1 \cup A_2$ has at most
    $(1 - \gamma - \zeta)n$ non-neighbors in $B_0$.
    Every $a \in A' \setminus (A_1 \cup A_2)$ has at most
    $2(1 - \gamma - \zeta)n$ non-neighbors in $B_0$.
  \end{claim}
  \begin{proof}
    Let $a \in A'$.
    First suppose $a \in A_i = N(c_i, A')$ for some $i \in \{1, 2\}$, then
    $a$ has no neighbors in $B_i$, so
    \begin{equation*}
      |N_{\overline{G}}(a, B_0)| \le
      |N_{\overline{G}}(a, B)| - |B_i| \le  
      (1 - \gamma - \zeta)n.
    \end{equation*}

    Now assume that $a \in A' \setminus (A_1 \cup A_2)$. 
    By Claim~\ref{clm:ac_ipath}, there exists $i \in \{1, 2\}$, 
    $c' \in C'$, and $a' \in A'$ such that $ac'a'c_i$ is a path.
    Because $ac'$ is an edge, $a$ has no neighbors in $N_G(c', B')$, thus 
    the number of non-neighbors of $a$ in $B\setminus N_G(c', B')$ is at most 
    $$|N_{\overline{G}}(a, B)| - |N_G(c', B')| \le |N_{\overline{G}}(a, B)|- \ceiling{(1/2 - y)n}\le
      (1 - \gamma - \zeta)n,$$    
    so the number of non-neighbors of $a$ in $B_0 \setminus N_G(c', B_0) \subseteq B \setminus N_G(c', B)$ 
    is at most $(1 - \gamma - \zeta)n$.
    To see that the number of non-neighbors of 
    $a$ in $N_G(c', B_0)$
    is at most $(1 - \gamma - \zeta)n$ (which proves the claim),
    note that 
    $N_G(c', B_0) \subseteq N_{\overline{G}}(a', B_0)$ 
    (because $G$ is triangle-free) and, by the first part of 
    the claim, the fact that $a' \in N(c_i, A') = A_i$ implies that
    $|N_{\overline{G}}(a', B_0)| \le (1 - \gamma - \zeta)n$.
  \end{proof}

  Now we will estimate $e(\overline{G}[A', B_0])$ from both sides. Recall that $A_1$ and $A_2$ are disjoint, so
  $|A_1 \cup A_2| \ge 2(\beta - x)n$.
  This with Claim~\ref{clm:path} implies
  \begin{equation}\label{eq:no3disjoint1}
    \begin{split}
      e(\overline{G}[A', B_0]) 
      &\le |A_1 \cup A_2|\cdot (1 - \gamma - \zeta)n + 
      |A' \setminus (A_1 \cup A_2)| \cdot 2(1 - \gamma - \zeta)n \\
      & \le 
      2(\beta - x)(1 - \gamma - \zeta)n^2 +
      (1 - 2\beta + x) \cdot 2(1 - \gamma - \zeta)n^2 \\
      &= 2(1 - \beta)(1 - \gamma - \zeta)n^2.
  \end{split}
  \end{equation}
  (In \eqref{eq:no3disjoint1}, we used that $1 - \gamma - \zeta \ge 0$, 
  which is implied by Claim~\ref{clm:path}.)
  By Claim~\ref{clm:2}, for every $b \in B_0$ there exists $c \in N(b, C')$. 
  Since $N_{\overline{G}}(b, A') \supseteq N(c, A')$, 
  \begin{equation}\label{eq:no3disjoint2}
    e(\overline{G}[A', B_0]) \ge |B_0|(\beta - x)n = 
    (1 - y - 2\zeta)(\beta - x)n^2.
  \end{equation}
  The conclusion then follows because 
  \eqref{eq:no3disjoint1} and \eqref{eq:no3disjoint2} together yield
  \begin{equation*}
    (1 - y - 2\zeta)(\beta - x) \le 2(1 - \beta)(1 - \gamma - \zeta)
  \end{equation*}
  which has no solutions when $x \ge 1/3$, $1/3 > y \ge \gamma - 1/2$, 
  $z \ge \beta - 1/4$, and $\zeta \ge 1/2 - y$.
  (See Lemma~\ref{lem:app5} in the appendix for a proof of this fact.)

\bibliographystyle{plain}

\appendix
\section{Proof of Lemma~\ref{lem:absorb}}
Let $V_1, \dotsc, V_k$ be the parts of $G$, let 
$\ell := k(t+1)$, and let
$\cA$ be the set of all $n^{\ell}$ sequences $a_1,\dotsc, a_{\ell}$ such that
$a_j \in V_i$ if $j$ is equivalent to $i$ modulo $k$.
Note that we do not require the vertices $a_1, \dotsc,a_{\ell}$ to be distinct in this definition,
so $|\cA| = n^\ell$.

For every transversal $U$,
define $\cA_U$ to be the set of sequences in $\cA$ such that 
if $A$ is the set of vertices in $\cA_U$, the graph induced by $A$ and 
the graph induced by $A \cup U$ both have a transversal $C_k$-factor.
The probabilistic argument below relies critically on
the fact that, for every transversal $U$, the set $\cA_U$ 
is sufficiently large, and this follows
from the fact that $G$ is $(2 \eta, t)$-linked.
To see this, 
first label the vertices in $U$ as $u_1, \dotsc, u_k$ so that
$u_i \in V_i$ for $i \in [k]$.
Because $G$ is $(2\eta,t)$-linked 
we easily have that  
there are at least $\eta n^k$ ways to select vertices $c_1, \dotsc, c_k$
where $c_i \in V_i \setminus \{u_i\}$ for $i \in [k]$
that induce a transversal $C_k$ in $G$.
Because $G$ is $(2\eta, t)$-linked, iteratively, for $i$ from $1$ to $k$, 
we can select a $(c_i, u_i, t)$-linking sequence $L_i$ 
that avoids all previously selected vertices in at least $\eta n^t$ ways.
Since the graph induced by $c_i$ and the vertices in $L_i$ contains a transversal $C_k$-factor
we have that $(t+1)$ is divisible by $k$ and 
that there exists $S_i$ an ordering of these $(t+1)$ vertices 
so that the $j$th vertex is in $V_i$ if $j$ is equivalent to $i$ modulo $k$.
Finally, because each $L_i$ is a $(c_i, u_i, t)$-linking sequence,
the concatenation of the sequences $S_1, \dotsc, S_k$
is in $\cA_U$, and we have that $|\cA_U| \ge \eta n^k \cdot (\eta n^t)^k = \eta^{k+1} n^{\ell}$. 

Let $p := 0.2 \cdot \sigma \cdot n^{-\ell + 1}$ 
and select the elements of $\cA$ independently with probability $p$
to form the random set $\cA_{\textrm{rand}}$.
The Chernoff and union bounds imply that with high-probability
\begin{equation}\label{eq:absorb_chernoff}
  |\cA_\textrm{rand}| \le \sigma n \qquad\text{and}\qquad 
  |\cA_\textrm{rand} \cap \cA_U| \ge 0.1 \cdot \sigma \eta^{k+1} n \ge (\ell^2 + 1)\sigma^2 n
\end{equation}
for every transversal $U \subseteq V(G)$.
Note that the number of pairs of sequences in $\cA$ in which a vertex is repeated is less than
$n \cdot \binom{2 \ell}{2} \cdot n^{2\ell - 2} = \binom{2 \ell}{2} n^{2\ell - 1}$, 
so the expected number of pairs of sequences in $\cA_\textrm{rand}$ 
in which a vertex is repeated is less than
$p^2 \cdot \binom{2 \ell}{2} n^{2\ell - 1} \le (\ell^2 \sigma^2 n)/4$. 
So, by Markov's inequality, with probability at least $1/2$, if we add both elements from
every such pair to form the set $\cA_\textrm{rep}$ we have that
\begin{equation}\label{eq:absorb_repetition}
  |\cA_\textrm{rep}| \le \ell^2 \sigma^2 n.
\end{equation}
Therefore, there exist an outcome in which both \eqref{eq:absorb_chernoff} and 
\eqref{eq:absorb_repetition} hold.
From $\cA_\textrm{rand}$ we now remove all sequences that are in $\cA_\textrm{rep}$
and all sequences for which there does not exists a transversal $U$ 
for which it is a linking sequence to form the collection $\cA'$.
Note that, by \eqref{eq:absorb_chernoff}, 
$z: = |\cA'| \le \sigma n$ and,
by \eqref{eq:absorb_chernoff} and \eqref{eq:absorb_repetition},
$|\cA' \cap \cA_U| \ge \sigma^2 n$
for every transversal $U \subseteq V(G)$. 
Let $A$ be the vertices that appear in a sequence of $\cA'$.
Because no vertex is repeated in $\cA'$ we have that
$|A \cap V_i| = z$ for every $i \in [k]$, 
and, because every sequences in $\cA'$ is a linking sequences for
some transversal $U$, for every sequence in $\cA'$ the graph
induced by the vertices in the sequence has a transversal $C_k$-factor.

Suppose that there exists a transversal $C_k$-tiling of $G - A$ that covers all of the 
vertices in $V(G - A)$ except a set $W$ such that $|W| \le k \sigma^2 n$.
We can arbitrarily partition $W$ into transversals $U_1, \dotsc, U_m$ where $m = |W|/k \le \sigma^2 n$.
Since for every $i \in [m]$, we have that $|\cA_{U_i} \cap \cA'| \ge \sigma^2 n \ge m$,
we can greedily select distinct sequences $A_1, \dotsc, A_m$
such that $A_i \in \cA_{U_i} \cap \cA'$ 
for every $i \in [m]$.
This implies that there is a transversal $C_k$-factor of $G[W \cup A]$ and, therefore,
a transversal $C_k$-factor of $G$.
\qedsymbol

\section{Inequalities from Section~\ref{sec:triangle_covers}}

\begin{lemma}\label{lem:app1}
  For every $1/2 \le \beta \le 2/3$ the following holds when 
  $\gamma = 4/3 - \beta$. 
  There does not exist $x,y,z$ such that $1 > x + y +z$, 
  $x \ge 1/3$, $y \ge 1/2$, $z \ge \beta - 1/2$, and
  \begin{equation*}
    (1 - \gamma)(z - \beta + 1/2) \ge (1/2 - z)(\beta - x).
  \end{equation*}
\end{lemma}
\begin{proof}
  Since $y \ge 1/2$ and $x \ge 1/3$ imply that
  $z < 1/6$ and $1/2 - z > 1/3 \ge 1 - \gamma$, we have 
  \begin{multline*}
    (1/2 - z)(\beta - x) > (1/2 - z)(\beta - (1 - y - z)) \ge \\
    (1/2 - z)(z + (\beta - 1/2)) \ge 
    (1/2 - z)(z - (\beta - 1/2)) > 
    (1 - \gamma)(z - (\beta - 1/2)). \qedhere
  \end{multline*}
\end{proof}

\begin{lemma}\label{lem:app2}
  For every $1/2 \le \beta \le 2/3$ the following holds when 
  $\gamma = 4/3 - \beta$. 
  There does not exist $x,y,z$ such that $1 > x + y +z$, 
  $x \ge 1/3$, $y \ge 5/12$, $z \ge \beta - 1/2$, and
  \begin{equation*}
    (1 - \gamma)(1 - z) \ge (1-x)(\beta - z).
  \end{equation*}
\end{lemma}

\begin{proof}
  Note that the conditions imply that $x+z< 7/12$ and $\gamma\ge 5/6-z$.
Assume there is a solution, so
  \begin{equation*}
  (1 - z) + \gamma(z-x)-(1-x)(4/3-z)=(1 - \gamma)(1 - z)-(1 - x)(4/3-\gamma - z)\ge 0.
  \end{equation*}
  Since $z-x<0$ (otherwise $x+y+z\ge2/3+5/12>1$), applying $\gamma\ge 5/6-z$ yields
  \begin{equation*}
       (1 - z) + (5/6 - z)(z - x) - (1 - x)(4/3 - z)\ge 0,
  \end{equation*}
 after simplification we have,
    \begin{equation*}
      (5/6 - z)z+x/2 - 1/3\ge 0,
  \end{equation*}
  using $x< 7/12-z$ we get,
        \begin{equation*}
      (1/3 - z)z - 1/24\ge 0
  \end{equation*}
  This has no solution, a contradiction. 
\end{proof}

\begin{lemma}\label{lem:app3}
  For every $1/2 \le \beta \le 2/3$ the following holds when 
  $\gamma = 4/3 - \beta$. 
  There does not exist $x,y,z$ such that $1 > x + y +z$, 
  $x \ge 1/3$, $y \ge 1/3$, $z \ge \beta - 1/2$, and following three inequalities hold
  \begin{equation}  \label{eq:lemma21-1}
    (1 - \gamma)(z - \beta + 1/2) \ge (\beta - z)(1/2 - y) 
  \end{equation}
    \begin{equation}  \label{eq:lemma21-2}
    (1 - \gamma)(z - \beta + 1/2) \ge (\beta - x)(1/2 - z)
  \end{equation}
    \begin{equation}  \label{eq:lemma21-3}
     (1 - \gamma)(1 - z) \ge (1-x)(\beta - z).
  \end{equation}  
  
  \end{lemma}
\begin{proof}

Assume for a contradiction that there is a solution.
We can then make following claims.
\begin{claim}
  \label{claim:zge16}
  $z>1/6$.
\end{claim}
  \begin{proof}
  Assume for a contradiction that $z\le 1/6$.
  Summing up \eqref{eq:lemma21-1} and \eqref{eq:lemma21-3} we get
    \begin{equation*}  
    (1 - \gamma)(3/2 - \beta ) \ge (\beta - z)(3/2 - x - y) 
  \end{equation*}
  Since $\beta=4/3-\gamma$, and $1-x-y>z$, we have
    \begin{equation*}  
    1/6 +5\gamma/6 - \gamma^2 > 2/3- \gamma/2 - z(z+ \gamma - 5/6)    
  \end{equation*}
  $z\ge \beta-1/2= 5/6-\gamma$, $z+ \gamma - 5/6\ge 0$, therefore $z(z+ \gamma - 5/6)$ is an increasing function with respect to $z$, so using $z\le 1/6$ we get
      \begin{equation*}  
    1/6 +5/6\gamma - \gamma^2> 7/9- 2\gamma/3 \implies
    -\gamma^2 + 3\gamma/2 - 11/18 > 0.
  \end{equation*}
 It is easy to check that this does not have a solution, a contradiction. 
  \end{proof}
  
\begin{claim}
  \label{claim:gammaplusz}
If $\gamma \le3/4$, then $\gamma+z>1$.
\end{claim}
  \begin{proof}
  $\beta-x=4/3-\gamma-x > 1/3-\gamma +y+z\ge 2/3-\gamma+z$, therefore by \eqref{eq:lemma21-2} 
  \begin{equation*} 
    (1 - \gamma)(z+\gamma -5/6) > (2/3-\gamma+z)(1/2 - z).
  \end{equation*}
  So,
    \begin{equation*} 
     z^2-\gamma^2-2\gamma z+7 \gamma/3 +7z/6-7/6>0.
  \end{equation*}
  The derivative of the above equation with respect to $z$ is
      \begin{equation*} 
    2z-2\gamma + 7/6 \ge 2\cdot 1/6 -2\cdot 3/4 + 7/6 = 0
  \end{equation*}
  Therefore, using $z<1-x-y\le 1/3$, we get
 \begin{equation*} 
    (1 - \gamma)(\gamma -1/2) - (1-\gamma)(1/2 - z)\ge (1 - \gamma)(z+\gamma -5/6) - (2/3-\gamma+z)(1/2 - z)> 0,
  \end{equation*}
 so $\gamma-1/2-1/2+z>0$.
  \end{proof}
  
Since $1-x>y+z\ge 1/3 + z$, \eqref{eq:lemma21-3} implies
\begin{equation*}
    (1 - \gamma)(1 - z) > (1/3+z)(4/3-\gamma- z).
\end{equation*}
by rearranging the terms we get
\begin{equation*}
(z-1/3)(2\gamma+z-5/3)>0.
\end{equation*}
Note that $z < 1 - x - y \le 1/3$,
so $z - 1/3 < 0$.
Therefore,
\begin{equation*}
2\gamma+z-5/3 < 0.
\end{equation*}
If $\gamma \le 3/4$, then 
Claim \ref{claim:gammaplusz} implies that $2\gamma+z-5/3>\gamma+1-5/3\ge 0$, a
contradiction. 
If $\gamma > 3/4$, then, with Claim~\ref{claim:zge16}, we have
$2\gamma+z-5/3 > 3/2 + 1/6 - 5/3 = 0$,
a contradiction.
%
%
%
%
%
\end{proof}

\begin{lemma}\label{lem:app4}
  For every $1/2 \le \beta \le 2/3$ the following holds when 
  $\gamma = 4/3 - \beta$. 
  There does not exist $x,y,z$ such that $1 > x + y +z$, 
  $x \ge 1/3$, $1/3 \ge y \ge \gamma - 1/2$, $z \ge \beta - 1/2$, and
  \begin{equation*}
    (1 - \beta)\left(y - \gamma + 1/2\right) 
    \ge
    (\gamma - y) \left(1/2 - z\right).
  \end{equation*}
\end{lemma}

\begin{proof}
Assume for a contradiction that there is a solution. $1/2-z> -1/2+x+y\ge y-1/6$, therefore
\begin{equation*}
    (1 - \beta)\left(y - \gamma + 1/2\right) 
    >
    (\gamma - y) \left(y-1/6\right).
  \end{equation*}
  Replacing $\beta$ by $4/3-\gamma$ and simplifying we get,
  \begin{equation*}
    y^2-y/2+\gamma-\gamma^2-1/6> 0.
  \end{equation*}

Note that $1/3\ge y \ge \gamma -1/2\ge 1/6$ implies $y^2-y/2\le -1/18$ and that $\gamma\ge 2/3$ implies $\gamma-\gamma^2\le 2/9$. Therefore, we get
$-1/18+2/9-1/6>0$, a contradiction.
\end{proof}

\begin{lemma}\label{lem:app5}
  For every $1/2 \le \beta \le 2/3$ the following holds when 
  $\gamma = 4/3 - \beta$. 
  There does not exist $x,y,z,\zeta$ such that $1 > x + y +z$, 
  $x \ge 1/3$, $1/3 \ge y \ge \gamma - 1/2$, $z \ge \beta - 1/4$,
  $\zeta \ge 1/2 - y$, and 
  \begin{equation*}
    (1 - y - 2\zeta)(\beta - x) \le 2(1 - \beta)(1 - \gamma - \zeta).
  \end{equation*}
\end{lemma}
\begin{proof}
Assume for a contradiction that the inequality holds.  Then,
  \begin{equation*}
    (1 - y)(\beta - x) \le 2(1 - \beta)(1 - \gamma) + \zeta(2\beta - 2x)-\zeta (2 - 2\beta).
  \end{equation*}  
  Which implies
  \begin{equation*}
    (1 - y)(\beta - x) \le 2(1 - \beta)(1 - \gamma) + \zeta(4\beta - 2x - 2).
  \end{equation*}  
  Since $\beta \le 2/3$ and $x \ge 1/3$, we have $4\beta - 2x - 2<0$.
  So, with the fact that $\zeta \ge 1/2 - y$, after rearranging we get
      \begin{equation*}
    y(\beta - x) \le 2(1 - \beta)(\beta - 1/3) + (1/2 - y)(2\beta - 2).
  \end{equation*}
  Plugging $x<1-y-z\le 5/4-y-\beta$ and simplifying, we have  
        \begin{equation*}
   11\beta/3 - 2\beta^2 + 13/4y-y^2 - 4\beta y - 5/3>0,
  \end{equation*} 
  which implies 
  \begin{equation*}
   11\beta/6 + 17y/12 + 11/6(\beta+y)-(\beta + y)^2-\beta^2 - 2\beta y - 5/3>0.
  \end{equation*}
  By the AM-GM inequality $(\beta+y)(11/6-(\beta +y))<121/144$, therefore
 \begin{equation*}
   11\beta/6 + 17y/12 + 121/144-\beta^2-2\beta y -5/3>0,
  \end{equation*}
  so
  \begin{equation*}
   11\beta/6 + y(17/12 - 2\beta) - \beta^2 - 119/144>0.
  \end{equation*}
 Since $17/12-2\beta >0$ and $y \le 1/3$,
 \begin{equation*}
   11\beta/6 + 17/36 - 2\beta/3 - \beta^2 - 119/144>0,
  \end{equation*}
  so
  \begin{equation*}
   7\beta/6 - \beta^2 - 51/144>0,
  \end{equation*}
   which does not have a solution, a contradiction.
  
\end{proof}
\end{document}